\documentclass[final,3p, 11pt]{elsarticle}

\usepackage{amsmath,amssymb,amsfonts,bbm,textcomp,amsthm,mathrsfs}

\usepackage{varioref}

\usepackage{mathtools}
\usepackage{comment}

\usepackage[h]{esvect}

\usepackage{array}

\usepackage{listings}
\usepackage[activate]{pdfcprot}
\usepackage{pdfpages}
\usepackage{booktabs}

\usepackage{textcomp}
\usepackage{tabulary}

\usepackage{multicol}
\usepackage{stmaryrd}

\usepackage{bbm}
\usepackage{bigints}
\usepackage{lmodern}

\DeclareFontFamily{U}{mathx}{\hyphenchar\font45}
\DeclareFontShape{U}{mathx}{m}{n}{
      <5> <6> <7> <8> <9> <10>
      <10.95> <12> <14.4> <17.28> <20.74> <24.88>
      mathx10
      }{}
\DeclareSymbolFont{mathx}{U}{mathx}{m}{n}
\DeclareFontSubstitution{U}{mathx}{m}{n}
\DeclareMathAccent{\widecheck}{0}{mathx}{"71}
\DeclareMathAccent{\wideparen}{0}{mathx}{"75}

\usepackage{hyperref}

\newcommand{\im}{\ensuremath{\mathrm{i}}}

\newcommand{\C}{\mathbb{C}}
\newcommand{\E}{{\mathbb E}}

\newcommand{\N}{\mathbb{N}}
\newcommand{\Pa}{{\mathbb P}}

\newcommand{\R}{\mathbb{R}}

\newcommand{\Acal}{{\mathcal A}}
\newcommand{\Bcal}{{\mathcal B}}
\newcommand{\Ccal}{{\mathcal C}}
\newcommand{\Dcal}{{\mathcal D}}

\newcommand{\Fcal}{{\mathcal F}}

\newcommand{\Hcal}{{\mathcal H}}

\newcommand{\Jcal}{{\mathcal J}}

\newcommand{\Lcal}{{\mathcal L}}
\newcommand{\Ucal}{{\mathcal U}}
\newcommand{\Mcal}{{\mathcal M}}

\newcommand{\Oh}{{\mathcal{O}}}

\newcommand{\Rcal}{{\mathcal R}}

\newcommand{\Rp}[1]{{\mathbb R^{#1}_{\geq0}}}
\newcommand{\Rpp}[1]{{\mathbb R^{#1}_{>0}}}

\newcommand{\expvBig}[2]{\operatorname{\E}^{#1}\Bigl[#2\Bigr]}

\newcommand{\tr}{\mathop{\mathrm{Tr}}}

\newcommand{\tc}{\widehat}

\newcommand{\la}{\left\langle}
\newcommand{\ra}{\right\rangle}

\newcommand{\st}{\;|\;}

\renewcommand{\epsilon}{\ensuremath\varepsilon}

\renewcommand{\phi}{\ensuremath{\varphi}}

\SetSymbolFont{stmry}{bold}{U}{stmry}{m}{n}

\newtheorem{theorem}{Theorem}[section]

\newtheorem{remark}[theorem]{Remark}
\newtheorem{corollary}[theorem]{Corollary}
\newtheorem{definition}[theorem]{Definition}
\newtheorem{lemma}[theorem]{Lemma}
\newtheorem{proposition}[theorem]{Proposition}
\newtheorem{notation}[theorem]{Notation}

\usepackage{lineno}

\journal{no name}

\begin{document}

\begin{frontmatter}



\title{Regularity Results for Degenerate Kolmogorov Equations of Affine Type }


\author{Nicoletta Gabrielli}

\address{Z\"urich, Switzerland}

\begin{abstract}
Motivated by option pricing problems, we are interested in the approximation of the quantity $\expvBig{x}{f(X_t)}$ where $X=(X_t)_{t\geq0}$ is a Markov process with $X_0=x$ 
and $f$ is a function possibly growing at 
infinity. Under a set of reasonable assumptions of compatibility between the vector field and the growth of the payoff function we can interpret the quantity as the 
action of a Markov semigroup which in turn can be identifies with the solution of the 
Kolmogorov PIDE driven by the (extended) generator of the process $X.$

\end{abstract}


\end{frontmatter}


\setlength{\parindent}{0pt}
\section{Introduction}

This paper is devoted to the study of Markov semigroups acting on spaces which consist of functions which are not necessarily 
bounded. More precisely, given a Markov process taking values in $D\subseteq\R^d$
$$(\Omega, (X_t)_{t\geq0},(\Fcal_t)_{t\geq0},(p_t)_{t\geq0},(\Pa^x)_{x\in D})$$
and a function $f\in \Mcal$ such that $\expvBig{x}{|f(X_t)|}<\infty$ for all $t\geq0$ and $x\in D,$ we want to find a set of 
conditions on the function space $\Mcal$ such that it is possible to identify the transition semigroup $P_t$ acting on $f$ with the solution of the Kolmogorov equation 
\begin{equation}\label{eq:ko}\tag{$\ast$}
 \begin{aligned}
   \partial_t u(t,x)&=\Acal u(t,x),&&\qquad (t,x)\in[0,T]\times D,\\
u(0,x)&=f(x),&&\qquad x\in D,
 \end{aligned}
\end{equation}
where $\Acal$ here denotes the (extended) generator of the Markov process $X.$

This connection between the Kolmogorov equation and the extended generator of the Markov process has been already derived in case 
the function space $\Mcal$ is given by a so called $\Bcal_{\psi}$ space for the Markov semigroup.  We refer to \cite{higa, doersek_semigroup_2010}, and 
\cite{doersek_semigroup_2010} for the definition of $\Bcal_\psi$ spaces and the theory of a generalized 
version of Feller property which allows unbounded payoff functions.

Unfortunately for affine processes 
there is not, at least at the moment, a suitable $\Bcal_\psi$ formulation of the problem and therefore the aforementioned approaches cannot be directly applied.

In this paper we focus on a martingale approach. Starting with the set of functions 
\begin{eqnarray*}
\Mcal&:=&\Big\{ f:D\to\R\st\mbox{ Borel measurable such that}\\
&& P_t |f|<\infty \mbox{ for all }t\geq0\mbox{ and }M^f\mbox{ is a true martingale}\Big\},
\end{eqnarray*}
where
\begin{equation*}
M^f_t:=f(X_t)-f(x)-\int_0^t \Acal f(X_s)ds,
\end{equation*} 
we show that, under the assumption that $\lim_{t\to0}P_t f(x)=f(x)$ for all $f\in\Mcal,$ it is possible to conclude that $P_t f(x)$ coincides with the solution 
$u(t,x)$ of the Kolmogorov equation $(\ast).$ 

In particular, we can conclude that, if $\Mcal$ is a set containing functions $f$ such that 
\begin{description}
\item[A1)] $P_t |f|<\infty,$
\item[A2)] $M^f$ is a true martingale,
\item[A3)] $t\mapsto P_t f(x)$ is continuous at $t=0,$
\end{description}
then $P_t\Mcal\subseteq\Mcal$ and it is possible to derive a Taylor expansion of the function $t\mapsto P_tf(x)$  for $t$ around zero.

However, in applications, we are given a set of functions $\Hcal$ and 
we search for conditions under which 
$P_t\Hcal\subseteq\Hcal.$ We first focus on the set $\Hcal$ of smooth functions with growth controlled by a weight function $F$ which satisfies
\begin{description}
\item[B1)] $P_t |F|<\infty,$
\item[B2)] $|\Acal F|\leq KF$ for some constant $K>0$ and $\expvBig{x}{\sup_{t\in[0,T]}F^2(X_t)}<\infty.$
\end{description}
Under these two conditions we will see that all the three conditions {\bfseries A1)}, {\bfseries A2)} and {\bfseries A3)} hold and therefore $P_t \Hcal\subseteq \Mcal$ and 
the Taylor expansion is possible. Actually condition {\bfseries B2)} together with Gronwall's lemma also implies $P_t\Hcal\subseteq \Hcal.$ However, 
the condition $\Acal\Hcal\subseteq\Hcal$ is not always possible to achieve, in particular when the definition of $\Hcal$ comprises 
also some type of order of differentiability for the test functions.

This considerations will lead to a set of conditions under which $t\mapsto P_t f(x)$ is differentiable up to a fixed order. Then, we would like to obtain an analogous result 
for $x\mapsto P_t f(x)$. We will see that, for affine processes, there exists a time--space transformation which allows us to translate any result about 
regularity with respect to the space variable in a statement about regularity in time.

The paper is organized as follows.

In Section \ref{sec:generalregu} we consider functions $f$ which belong to the domain of the extended generator and provide a framework for the 
Kolmogorov equation to make sense. 

In Section  \ref{sec:weighted} we restrict ourselves to the space $\Hcal$ of all the functions which are infinitely differentiable and with growth 
controlled by a weight function. 
We analyze under which conditions the function $u:\Rp{}\times D\to\R$ defined by 
$$u(t,x):=P_t f(x)=\expvBig{x}{f(X_t)}$$ has derivatives of all orders satisfying the following property
$$\mbox{for all }(t,x)\in[0,T]\times D,\quad \partial^\alpha_{(t,x)}u(t,x)\in\Hcal.$$ 

In Section 4 we will also exploit in details the conditions we obtain when the Markov process $X$ is a L\'evy process. 

We conclude with some applications to weak approximation of the trajectories of an affine process. 

In the last section we apply the results derived in Section \ref{sec:generalregu} and Section \ref{sec:weighted} to the class of 
affine processes and the function space $\Hcal=C^\infty_{\textrm{pol}}$ of smooth functions with polynomial growth.

\paragraph{Notation}

Henceforth $D$ denotes the subset $\Rp{m}\times\R^n$ of $\R^d.$ The scalar 
product is denoted by $\la\cdot,\cdot\ra$ and associated norm $|\cdot|.$ 
The same notation is used also when the scalar product is considered in the space $\R^d+i \R^d.$ In this case we mean the extension of $\la\cdot,\cdot\ra$ 
in $\R^d+i\R^d$ without conjugation.

Given a measure $\mu$ taking values on $D$ the quantity 
\begin{equation*}
\int_D e^{\la u,\xi\ra}\mu(d\xi)
\end{equation*}
 is well defined for all 
$x\in D, t\geq0$ if and only if $u$ lies in the set 
\begin{equation}\label{ucal}
 \Ucal:=\left\{u\in\C^d\st x\mapsto e^{\la x,u\ra}\mbox{ is bounded on }D\right\}.
\end{equation}

Due to the geometry of the state space, the function
\begin{equation}\label{eq:fux}
f_u(x):=e^{\la x,u\ra}
\end{equation}
is bounded if and only if $\Ucal:=\C^m_{\leq0}\times \im\R^n.$

\begin{notation}\label{not:proj}
 In order to simplify the notation we introduce the following sets of indices $I$ and $J$ 
defined as $I=\{1,\ldots,m\}$ and $J=\{m+1,\ldots,m+n\}.$
Moreover, given a set $H\subseteq\{1,\ldots,d\}$ the map $\pi_H$ is the projection of the canonical state space on the lower dimensional subspace indexed by $H.$ 
In particular
\begin{eqnarray*}
{\pi_I}:\Rp{m}\times\R^n&\to&\Rp{m}\\
x&\mapsto&\pi_I x:=(x_i)_{i\in I}\\
\end{eqnarray*}
\begin{eqnarray*}
{\pi_J}:\Rp{m}\times\R^n&\to&\R^n\\
x&\mapsto&\pi_J x:=(x_j)_{j\in J}.
\end{eqnarray*}
\end{notation}

\begin{definition}\label{def:infdivD}
A function $\eta:\Ucal\to\C$ has the L\'evy-Khintchine form on $D=\Rp{m}\times\R^n$ 
if there exists $\mathsf{b}\in D, \sigma\in S^n_+$ and a 
Borel measure $\nu$ on $D$ satisfying 
$$\nu(\{0\})=0\quad\mbox{ and }\quad \int_{D}((|\pi_I\xi|+|\pi_J\xi|^2)\wedge 1)\nu(d\xi)<\infty$$ 
such that for any $u\in\C^m_{\leq0}\times\im\R^n$
\begin{align}
\eta(u)&=\la \mathsf{b},u\ra+\frac{1}{2}\la \pi_J u,\sigma \pi_J u\ra\nonumber\\
&+\int_{D\setminus\{0\}}\left(e^{\la u,\xi\ra}-1-\la \pi_J u,\pi_J h(\xi)\ra\right)\nu(d\xi),\label{eq:infdiv}
\end{align}
where  $\pi_I$ and $\pi_J$ denote the projections of elements of $\Rp{m}\times\R^n$ into $\Rp{m}$ and $\R^n$ 
respectively, see Notation \ref{not:proj}.

\end{definition}

It is well known that there exists a unique infinitely divisible distribution in $D$ such that  
$$\int_D e^{\la u,\xi\ra}\mu(d\xi)=e^{\eta(u)},
$$
if and only if $\eta$ is of type \eqref{eq:infdiv}.

Let $\eta$ of type \eqref{eq:infdiv} and $f_u$ defined in \eqref{eq:fux}. Denote by $L^\eta$ the L\'evy process with L\'evy exponent $\eta,$ i.e.
\begin{equation}\label{eq:deflevy}
\expvBig{}{f_u(L^\eta_s)}=e^{s\eta(u)},\qquad s\geq0,\; u\in\Ucal.
\end{equation}

We collect now the assumptions and conventions we make on Markov processes for this paper. 

Given a probability space $\Omega,$ let 
$$X=(\Omega,(X_t)_{t\geq0},(\Fcal_t)_{t\geq0},(p_t)_{t\geq0},(\Pa^x)_{x\in D}),$$ 
be a time-homogeneous Markov process with state space $D$ where
\begin{itemize}
\item $(X_t)_{t\geq0}$ stochastic process taking values in $D$
\item $\Fcal_t=\sigma(\{X_s\;, s\leq t\}),$ 
\item $(p_t)_{t\geq0}$ semigroup of transition functions on $(D,\Bcal(D)),$
\item $(\Pa^x)_{x\in D}$ probability measures on $(\Omega,\Fcal).$ 
\end{itemize}
satisfying
\begin{equation}\label{eq:markovprop}
\expvBig{x}{f(X_{t+s})\big|\Fcal_t}=\expvBig{X_t}{f(X_s)},\quad \Pa^x\mbox{-a.s. for all }f\in mbdd(D),
\end{equation}
where $mbdd(D)$ is the space all all the functions $f:D\mapsto \R$ which are measurable and bounded.

Given $f\in mbdd(D)$  denote by 
$$P_t f(x):=\expvBig{x}{f(X_t)}=\int_D f(\xi)p_t(x,d\xi)$$
the transition semigroup of the Markov process $X.$

Whenever we have a c\`adl\`ag Markov process, we will always consider the canonical version realized on the filtered space. 
In this case the law of $X$ under
 $\Pa^x$ is a probability measure on $\Dcal(D)$ and therefore we can assume, without loss of generality, that $\Pa^x$ is a 
measure on $\Dcal(D).$

Let $\Mcal$ be the set of all Borel measurable functions 
$f:D\mapsto \R$ such that the following integral is well defined
$$P_t |f|(x):=\int_D |f(\xi)|p_t(x,d\xi),\qquad x\in D,\,t\geq0.$$
We assume that the family of measures $\Pa^x$ is such that, for all $f\in\Mcal$, it holds
$$\expvBig{x}{f(X_{t+s})\Big|\Fcal_s}=\expvBig{X_s}{f(X_t)}=P_t f(X_s),\qquad \Pa^x\mbox{-a.s.}$$
for all $x\in D$ and $s,t\in\Rp{}.$

\begin{notation}
 Depending on the situation, we will use different notations to denote the same quantity
\begin{equation}\label{eq:euP}
\expvBig{x}{f(X_t)}=u(t,x)=P_tf(x).
\end{equation}
In particular, $u$ will be used when the function $f$ is fixed and we are interested in the analysis of the function 
$\expvBig{x}{f(X_t)}$ as a function of both variables $(t,x)\in\Rp{}\times D.$ The notation in 
terms of the transition semigroup will be particularly handy to compactify the notation for the time evolution dynamics. 
\end{notation}

\section{Martingale problem and short--time asymptotic formula}\label{sec:generalregu}
In this section we fix the problem of giving a comprehensive framework for the transition semigroup of a Markov process acting on a class of functions which is as general as possible. We follow the methodology used in \cite{polypro}.
We start with the definition of the 
\emph{extended generator}.
\begin{definition}[Definition 7.1. in \cite{cinlar_semimartingales_1980}]\label{def:extgen}
 Given a Markov process $X$, the \emph{extended generator} is an operator $\Acal$ with domain 
$\Dcal(\Acal)$ such that, for any $f\in\Dcal(\Acal)$ the process 
\begin{equation}\label{eq:Mf}
M^f_t:=f(X_t)-f(x)-\int_0^t \Acal f(X_s)ds
\end{equation}
is a local martingale under $\Pa^x$ for every $x\in D.$
\end{definition}

Consider a function $f\in\Dcal(\Acal)$ such that $\expvBig{x}{|f(X_t)|}<\infty$, for all $(t,x)\in\Rp{}\times D$ and such that $(M^f_t)_{t\geq0}$ is a $\Pa^x$ true 
martingale for every $x\in D$.  Then we can take the expectation on both side of \eqref{eq:Mf} 
$$\expvBig{x}{f(X_t)-f(x)-\int_0^t \Acal f(X_s)ds}=0\, .$$
Write
\begin{equation}
\begin{aligned}
 \expvBig{x}{f(X_t)}&=f(x)+\expvBig{x}{\int_0^t \Acal f(X_s)ds}\,,
\end{aligned}
\end{equation}
and define $P_t f(x):=\expvBig{x}{f(X_t)}$ so that the previous equation reads
$$P_t f(x)=f(x)+\int_0^t P_s \Acal f(x)ds\, .$$

From Fubini's theorem applied to the increments of the martingale $M^f$ 
(see Remark 4.1.4) in \cite{polypro}, it follows that 
 $\expvBig{x}{|\Acal f(X_s)|}<\infty$, for all $x\in D$ and $s\geq0$,  
 and therefore we have got an integral equation for $P_tf$. Hereafter, we restrict ourselves on $\Mcal$, defined as the space of functions 
$f$ such that $M^f$ is a true martingale:
\begin{equation}\label{eq:Mcal}
\begin{aligned}\Mcal:=\Big\{& f:D\to\R\st\mbox{ Borel measurable such that}\nonumber\\
& P_t |f|<\infty \mbox{ for all }t\geq0\mbox{ and }M^f\mbox{ is a true martingale}\Big\}.
\end{aligned}
\end{equation}

However, in order to conclude that $P_s\Acal f(x)=\Acal P_s f(x)$
for all $x\in D$ and $s\geq0$,  we need some additional regularity in time for the process.

In the next theorem we collect 
some results from 
\cite{polypro}.

\begin{theorem}[Lemma 2.6. in \cite{polypro}]\label{teo:FK} Let $X$ be a time-homogeneous Markov 
process and $f\in\Mcal$. Then 
\begin{enumerate}[(i)]
\item  For any $s\geq0$,  $M^{P_s f}$ is a true martingale.
\item If $t\to P_t f$ is continuous at $t=0$,   $\Acal P_t f=P_t \Acal f$ for any $t\geq0$.
\item (Feynman--Kac representation) If $t\to P_t\Acal f(x)$ is continuous at $t=0$,  then $P_t f(x)$ 
coincides with the solution $u(t,x)$ of the Kolmogorov equation
\begin{displaymath}
 \begin{array}{rcll}
   \partial_t u(t,x)&=&\Acal u(t,x),&\quad (t,x)\in[0,T]\times D\,,\\
u(0,x)&=&f(x),&\quad x\in D\,.
 \end{array}
\end{displaymath}

\end{enumerate}
\proof
This is simply an adaptation of the proof in \cite{polypro}. Since
$$f(X_t)-f(x)-\int_0^t \Acal f(X_s)ds$$
is a true martingale, all its increments  have vanishing expectation. Hence both $f(X_t)$ and $\Acal f(X_t)$ are
 integrable for every $t\geq0$. Consider the process
$$\widetilde M^{P_s f}_t:=P_s f(X_t)-P_s f(x)-\int_0^t P_s\Acal f(X_r)dr\, .$$
Due to the integrability of the increments of $M^f$, we have that $\widetilde M^{P_s f}$ is well 
defined for all $s$ and for all $t\in[0,T]$. Using Markov property of $X$ it holds, for $t_0\leq t$
\begin{align*}
 \expvBig{x}{\widetilde M^{P_s f}_t-\widetilde M^{P_s f}_{t_0}|\Fcal_{t_0}}=&\expvBig{x}{P_sf(X_t)-P_s f(X_{t_0})-\int_{t_0}^t P_s\Acal f(X_r)dr\big|\Fcal_{t_0}}\\
=&\expvBig{X^x_{t_0}}{P_s f(X_{t-t_0})-P_s f(X_{t_0})-\int_{t_0}^t P_s \Acal f(X_{r-t_0})dr}\\
=&\left(P_{s+t-t_0 }f(y)-P_s f(y)-\int_{s}^{s+t-t_0}P_r\Acal f(y)dy\right)_{\Big| y=X^x_{t_0}},
\end{align*}
but the last term is identically zero from martingale property of $M^f.$
Taking $t_0=0$,  we get that the map $r\mapsto P_s \Acal f(X_r)$ is integrable for all $s\in[0,t]$ for all $t$ and 
$$\int_0^t P_s|\Acal f(X_r)|dr<\infty,\quad\mbox{for all }t\in[0,T]\, .$$ By definition of extended generator, together with the fact that 
$P_sf$ is integrable, we conclude that $P_s f\in\Dcal(\Acal)$ and 
$$P_s\Acal f=\Acal P_s f\, .$$ Finally, this also implies that $M^{P_s f}$ is a true martingale for all $s\geq0$. Kolmogorov's equation follows by taking the limit of the finite differences and using continuity of $t\mapsto P_t\Acal f(x)$. Precisely
\begin{eqnarray*}
\partial_t P_tf(x)&=&\lim_{h\to0}\frac{P_{t+h}f(x)-P_t f(x)}{h}\\
&=&\lim_{h\to0}P_t\left(\frac{P_h f(x)-f(x)}{h}\right)\\
&=&\lim_{h\to0}P_t\frac{1}{h}\int_0^h P_s\Acal f(x)ds\\
&=&P_t\Acal f(x)=\Acal P_t f(x)\,,
\end{eqnarray*}
where in the last line we used commutative property between $P$ and $\Acal$. 
\endproof
\end{theorem}

Hence, under continuity assumption, both the formulations
$$P_t f(x)=f(x)+\int_0^t P_s \Acal f(x)ds$$
and
$$P_t f(x)=f(x)+\int_0^t \Acal P_s f(x)ds$$
are well defined for all $t\geq0$. In line with the literature, we refer to the first one as Dynkin's formula and the second one as Kolmogorov equation.

Combining together martingale property and commutative property between the transition semigroup 
and the extended generator we get
\begin{corollary}[Dynkin's formula]
Under the assumptions in Theorem \ref{teo:FK}, for any $t,s\geq0$, it holds
\begin{equation}\label{eq:dynkin}
\expvBig{x}{f(X_t)}=f(x)+t\int_0^1 \expvBig{x}{\Acal f(X_{rt})}dr.
\end{equation}
\end{corollary}

Dynkin's formula can be iterated by taking the successive powers of the infinitesimal generator. 
Define iteratively:
\begin{equation}\label{eq:iterop}
\begin{aligned}
 \Acal^0f&=f\,,\\
\Acal^{n+1}f&=\Acal(\Acal^n f)\quad\mbox{for }n\geq0\,.
\end{aligned}
\end{equation}
Then we following holds

\begin{proposition}[Iterated Dynkin's formula]\label{prop:dynkiniter}
 Let $X$ be a Markov process and $\nu\in\N$. For all $f\in\Dcal(\Acal)$ such that, for all $n=0,\ldots,\nu$,  
 $\Acal^{n+1}f\in\Mcal$,  it holds
\begin{equation}\label{eq:dynkiniter}
\begin{aligned}
 \expvBig{x}{f(X_t)}=&f(x)+\sum_{k=1}^\nu \frac{t^k}{k!}\Acal^k f(x)\\
 &+\frac{t^{\nu+1}}{\nu!}\int_0^1(1-r)^\nu \expvBig{x}{\Acal^{\nu+1}f(X_{rt})}dr\,.
 \end{aligned}
\end{equation}

\proof
The proof is done by induction on $\nu$. For $\nu=0$ 
$$\expvBig{x}{f(X_t)}=f(x)+t\int_0^1\expvBig{x}{\Acal^{1}f(X_{st})}ds,
$$
coincides with \eqref{eq:dynkin}. Suppose now that the formula holds for $\nu-1$ and we prove it for $\nu.$
We can write the right hand side in \eqref{eq:dynkiniter} as
\begin{align*}
& f(x)+\sum_{k=1}^{\nu}\frac{t^k}{k!}\Acal^k f(x)+\frac{t^{\nu+1}}{\nu!}\int_0^1 (1-s)^\nu\expvBig{x}{\Acal^{\nu+1}f(X_{st})}ds\\
&= f(x)+\sum_{k=1}^{\nu-1}\frac{t^k}{k!}\Acal^k f(x)+\frac{t^\nu}{\nu!}\Acal^\nu f(x)+\frac{1}{\nu!}\int_0^t (t-s)^\nu P_{s}\Acal^{\nu+1}f(x)ds\\
&=f(x)+\sum_{k=1}^{\nu-1}\frac{t^k}{k!}\Acal^k f(x)+\frac{t^\nu}{\nu!}\Acal^\nu f(x)\\
&+\Big[\int_0^t \frac{(t-s)^{\nu-1}}{(\nu-1)!}P_s\Acal^\nu f(x)ds-\frac{t^\nu}{\nu!}\Acal^\nu f(x)\Big],
\end{align*}
where in the last step we did integration by parts on the integral term. 

\endproof
\end{proposition}

Clearly one could try to use the same argument starting from the Kolmogorov equation. The additional differentiability assumption in the next result is necessary in order to permute the extended generator 
with the transition semigroup
\begin{corollary}[Short--time asymptotic formula]
Under the same assumptions as in Proposition \ref{prop:dynkiniter}, if moreover $t\mapsto P_t\Acal^{1}f$ is continuous at 
$t=0$  the following expansion of 
the transition semigroup holds:
\begin{equation}
\begin{aligned}
P_tf(x)=&f(x)+\sum_{k=1}^\nu \frac{t^k}{k!}\Acal^k f(x)\\
&+\frac{t^{\nu+1}}{\nu!}\int_0^1 (1-s)^\nu\Acal^{\nu+1}P_{st}f(x)ds.
\end{aligned}
\end{equation}

\end{corollary}

\section{Markov semigroup on weighted spaces}\label{sec:weighted}

\subsection{Functions with controlled growth}

\begin{definition}
Given a left-continuous, increasing function $\rho:\Rp{}\mapsto\Rpp{}$ with 
$\lim_{u\to\infty}\rho(u)=+\infty$,  fix $\eta\in\Rp{}$ and define 
\begin{equation}\label{eq:wfn}
\begin{aligned}
 F_\eta(x):D&\to\Rp{}\\
x&\mapsto\rho(\eta|x|).
\end{aligned}
\end{equation}
A function $f:D\to\R$ is a continuous function with growth controlled by $F_\eta$ 
if there exists a constant $C$ such that
$$|f(x)|\leq CF_\eta(x)\, .$$
 The space of all continuous functions with growth controlled by $F$ will be denoted by $
C_{F}:$
\begin{equation}
 C_{F}:=\Big\{f\in C\st \exists C>0,\,\eta>0\, |f(x)|\leq CF_\eta(x),\,\mbox{ for all }x\in D\Big\}.
\end{equation}

\begin{definition}
 Observe that, for each $f\in C_F$,  there exists a couple $(C,\eta)$ such that $|f(x)|\leq CF_\eta(x)$. We call $(C,\eta)$ a \emph{good couple} for $f.$
\end{definition}

 In the 
space $C_{F}$ we introduce the norm 
$$||f||_{C_{F}}:=\inf\{C>0\st\;(C,\eta)\mbox{ is a good couple for }f\mbox{ and } |f(x)|\leq C F_\eta(x)\}\, .$$
\end{definition}

\begin{lemma}\label{lem:continzerogeneral}
 Let $X$ be a time homogeneous Markov process. 
 Given $f\in C_F$ with good couple $(||f||_{C_F},\eta)$,  suppose that, for all $t\in[0,T]$  and 
$x\in D$,  it holds $\expvBig{x}{F_\eta(X_t)}<\infty.$
Then
$$\lim_{t\to 0^+}\expvBig{x}{f(X_t)}=f(x)\, .$$
\proof
For any $x\in D$,  let $R$ be a constant such that $|x|<R$. 
We decompose 
\begin{eqnarray*}
 \big|\expvBig{x}{f(X_t)}-f(x)\big|&\leq&\expvBig{x}{|f(X_t)-f(x)|\mathbbm1_{\{|X_t|\leq, R\}}}\\
&&+\expvBig{x}{|f(X_t)|\mathbbm1_{\{|X_t|>R\}}}\\
&&+f(x)\Pa^x(|X_t|>R).
\end{eqnarray*}
The first term can be made arbitrarily small by weak convergence. 
For $R$ big enough it holds
$$\expvBig{x}{|f(X_t)|\mathbbm1_{\{|X_t|>R\}}}\leq ||f||_{C_{F}}\expvBig{x}{F_\eta(X_t)\mathbbm1_{\{|X_t|>R\}}}\, .$$
Moreover, 
$$\Pa^x(|X_t|>R)\leq\frac{1}{\rho(\eta R)}\expvBig{x}{F_\eta(|X_t|)}\, .$$
Both terms go to zero as $R$ goes to infinity.
\endproof
\end{lemma}

\subsection{Differentiable functions with controlled growth}
Recall that, given a multi--index $\alpha=(\alpha_1,\ldots,\alpha_d)\in\N^d$,  with the symbol $\partial^\alpha_x$ we denote the mixed partial 
derivative of order $|\alpha|=\alpha_1+\ldots+\alpha_d$
$$\partial^\alpha_x:=\frac{\partial^{|\alpha|}}{\partial^{\alpha_1}_{x_1}\partial^{\alpha_2}_{x_2}\cdots\partial^{\alpha_d}_{x_d}}\, .$$

\begin{definition}
Fix a weight function $F$ and a constant $\eta$. 
A function $f$ is $k$-times differentiable with growth controlled by 
$F$ if, for each multi--index $\alpha\in\N^d$ with $|\alpha|\leq k$ 
there exist two constants $C_\alpha>0$ and $\eta_\alpha>0$ with $\eta_\alpha\in(0,\eta)$ 
such that, for all $x\in D,$
\begin{equation}\label{eq:boundder}
|\partial^\alpha_x f(x)|\leq C_\alpha F_{\eta_\alpha}(x).
\end{equation}
The space of all the functions which are $k$-times differentiable with growth controlled by 
$F$ will be denoted by $C^k_F:$
\begin{equation}
 \begin{aligned}
  C^k_F:=&\Big\{ f\in C^k\st \mbox{ for all }\alpha\in\N^d\mbox{ with } |\alpha|\leq k\,,\\
&\quad\exists\, C_\alpha>0, \eta_\alpha>0\mbox{ such that } |f(x)|\leq C_\alpha F_{\eta_\alpha}(x)\;\mbox{ for all }x\in D\Big\}\,.
 \end{aligned}
\end{equation}

 Therefore, for any $\alpha\in\N^d$, there exists a couple 
$(C_\alpha,\eta_\alpha)$ satisfying \eqref{eq:boundder}. Let $\eta_\infty:=\max_\alpha \eta_\alpha$. The following norm
$$||f||_{C^k_{F}}:=\inf\{C>0\st |\partial^\alpha_x f(x)|\leq C F_{\eta_\infty}(x)\quad\mbox{for all }\alpha\in\R^d\mbox{ with }|\alpha|\leq k\}$$
is well defined. This definition extends for smooth functions. We define 
\begin{equation}
\begin{aligned}
C^\infty_F&:=\{f\in C^\infty\st\mbox{ for all }\alpha\in\N^d\,, \\
 &\quad\exists\, C_\alpha>0, \eta_\alpha>0\mbox{ such that } |\partial^\alpha _x f(x)|\leq C_{\alpha} F_{\eta_\alpha}(x)\mbox{ for all }x\in D\}\,.
\end{aligned}
\end{equation}
In line with \cite{alfonsi_high_2010}, we call $(C_\alpha,\eta_\alpha)_{\alpha\in\N^d}$ a good sequence. Given $f\in C^\infty_F$ 
with good sequence $(C_\alpha,\eta_\alpha)$,  necessarily there exists $\eta_\infty:=\max_\alpha\eta_\alpha$. This implies, in 
particular, that it is possible to find a constant $\eta_\infty$ such that all the partial derivatives of any 
orders have growth controlled by $F_{\eta_\infty}$. Therefore the following norm is well defined
$$||f||_{C^\infty_{F}}:=\inf\{C>0\st |\partial^\alpha f(x)|\leq C F_{\eta_\infty}(x)\quad\mbox{for all }\alpha\in\R^d\}\, .$$

\end{definition}

Given a function $f\in C^\infty_F$,  we want to exploit conditions under which the process $M^f$ defined in \eqref{eq:Mf} 
is actually a true martingale. We first add conditions under which, given a weight function $F$, it holds $C^\infty_F\subseteq\Dcal(\Acal)$. Then, 
by definition, $M^f$ is a local martingale. Obviously, some additional conditions need to be added in order to conclude 
that the process is a true martingale. We start here with a set of some general conditions which guarantee 
square integrability of $M^f$. 

\begin{proposition}\label{prop:checkmg}
Let $X$ be a time--homogeneous Markov process with extended generator $\Acal$. Suppose that for some $\eta^*>0$
\begin{enumerate}[(i)]
\item there exists a constant $K$ such that, for all $x\in D,$
$$|\Acal f(x)|\leq KF_{\eta^*}(x)\, ,$$
\item for all $x\in D$, it holds
$$\expvBig{x}{\sup_{t\in[0,T]}F^2_{\eta^*}(X_t)}<\infty\, .$$
Then, for all $f\in C^\infty_F$ such that  $\eta_\infty <\eta^*$,  the process 
$M^f$ is a true martingale.
\end{enumerate}
  \proof
Since $f\in\Dcal(\Acal)$ and $P_t|f|<\infty$,  by definition of extended generator, $M^f$ is a local martingale. 
Hence, there exists an increasing sequence of stopping times with $\lim_{n\to\infty}\tau_n=\infty\,$ $\Pa^x-$a.s.  such that 
$(M^f_{t\wedge\tau_n})_{t\geq0}$ are martingales for all $n\in\N$. Since $f\in C^\infty_F$,  there exist two constants 
$C>0$ and $\eta_\infty>0$ such that 
$|f(x)|\leq C F_{\eta_\infty}(x)$. Henceforth, $\widetilde C$ is a constant which may vary from line to line.
For $t\in[0,T]$, it holds
\begin{align*}
|M^f_{t\wedge\tau_n}|^2=&\left| f(X_{t\wedge\tau_n})-f(x)-\int_0^{t\wedge\tau_n}\Acal f(X_s)ds\right|^2\\
\leq&\widetilde C \left(F^2_{\eta_\infty}(X_{t\wedge\tau_n})+F^2_{\eta_\infty}(x)+\int_0^{t}F^2_{\eta^*}(X_{s\wedge\tau_n})ds\right).
\end{align*}
Taking the expectations
\begin{align*}
 \expvBig{x}{|M^f_{t\wedge\tau_n}|^2}\leq&\widetilde C\left(\expvBig{x}{F^2_{\eta_\infty}(X_{t\wedge \tau_n})}+\expvBig{x}{F^2_{\eta_\infty}(x)}+\expvBig{x}{\int_0^t F^2_{\eta^*}(X_{s\wedge\tau_n})ds}\right)\\
\leq&\widetilde C\left(1+t\expvBig{x}{\sup_{s\in[0,t]} F^2_{\eta^*}(X_s)}\right).
\end{align*}
Using the second assumption, we see that, for all $x\in D$,  there exists a constant $C_x$ such that, for all $n\in \N$ and $t\in[0,T]$,
$$\expvBig{x}{|M^t_{t\wedge\tau_n}|^2}\leq C_x\, .$$
Using Doob's inequality we conclude that 
$$\expvBig{x}{\sup_{t\in[0,T]}|M^f_{t\wedge\tau_n}|}\leq \widetilde C\expvBig{x}{|M^f_{T\wedge\tau_n}|^2}\, .$$
By monotone convergence theorem we conclude that $$\expvBig{x}{\sup_{t\in[0,T]}|M^f_{t}|^2}<\infty\, ,$$
from where square integrability of the $M^f$ follows.
\endproof

\end{proposition}

\subsection{L\'evy--type operators on weighted spaces}

Let $L$ be a L\'evy process with L\'evy triplet  $(\mathtt b,\sigma,\nu)$. Henceforth, the extended generator of a L\'evy process is denoted by $\Lcal.$

For all  $f\in\Dcal(\Lcal),$
\begin{align*}
 \Lcal f(x)&=\la \mathtt b, \nabla f(x)\ra+\frac{1}{2}\tr(\sigma D^2 f(x))\\
&+\int_{D\setminus\{0\}}\left( f(x+\xi)-f(x)-\la h(\xi),\nabla f(x)\ra \right)\nu(d\xi),
\end{align*}

where $h(\xi)$ is a fixed truncation function. 
\begin{proposition}\label{prop:levydynk}
Let $L$ be a L\'evy process on $\R^d$ with L\'evy triplet $(\mathtt b,\sigma, \nu)$ and denote by $\Lcal$ its extended generator.  Assume that there exists $\eta^*$ such that 
$$\int_{\{|\xi|\geq1\}}(|\xi|^2\wedge F_{\eta^*}(\xi))\nu(d\xi)<\infty\, .$$ 
For all $f\in C^\infty_F$ such that $\eta_\infty<\eta^*$
\begin{enumerate}[(i)]
 \item  $f\in\Dcal(\Lcal)$ and $\Lcal f\in C^{\infty}_F,$
\item if $M^f$ is a true martingale, for any fixed $\nu\in\N$,  the following expansion holds
\begin{equation}\label{eq:taylorlevy}
\begin{aligned}
 \expvBig{y}{f(L_s)}&=f(y)+\sum_{n=1}^\nu \frac{s^n}{n!}\Lcal^h f(y)\\
 &+\frac{s^{\nu+1}}{\nu!}\int_0^1(1-r)^\nu \expvBig{y}{\Lcal^{\nu+1}f(L_{rs})}dr,
 \end{aligned}
\end{equation}
where the operators $\Lcal^n,\;n=1,\ldots,\nu+1$ are defined in \eqref{eq:iterop}.
\end{enumerate}
\proof

From Theorem I.4.57 in \cite{jacod_limit_1987}, $f(L)$ reads
\begin{eqnarray*}
 f(L_s)&=&f(y)+\int_0^s \la\nabla f(L_{r^-}),\mathtt b\ra dr+\int_0^s\sum_{i=1}^d\partial_{y_i}f(L_{r^-})\sum_{j=1}^d \sqrt{\sigma_{ij}}dW^j_r \\
&&+\frac{1}{2}\int_0^s \tr(D^2 f(L_{r^-})\sigma)dr\\
&&+\int_0^s \int\left( f(L_{r^-}+\xi)-f(L_{r^-})\right)\left(\Jcal^L(dr,d\xi)-\nu(d\xi)dr\right)\\
&&+\int_0^s \int \left(f(L_{r^-}+\xi)-f(X_{r^-})-\la h(\xi),\nabla f(L_{r^-})\ra\right)\nu(d\xi)dr\,.
\end{eqnarray*}
Since the second and the forth term on the right hand side are predictable processes of finite variation, 
$f(L)$ is a special semimartingale.
Hence
\begin{align*}
M^f_s&=f(y)+\int_0^s\sum_{i=1}^d\partial_{y_i}f(L_{s^-})\sum_{j=1}^d (\sqrt{\sigma})_{ij}dW^j_r\\
&+\int_0^s \int\left( f(L_{r^-}+\xi)-f(L_{r^-})\right)(\Jcal^L(dr,d\xi)-\nu(d\xi)dr)\,,
\end{align*}
is a true martingale by assumption and the extended generator $\Lcal$ is given by
 \begin{align*}
 \Lcal f(y)&=\la \mathtt b,\nabla f(y)\ra+\frac{1}{2}\tr(\sigma D^2 f(y))\\
&+\int\left( f(y+\xi)-f(y)-\la h(\xi),\nabla f(y)\ra\right)\nu(d\xi).
\end{align*}
The following estimate holds
$$||\Lcal f||_{C^{\infty}_{F}}\leq K||f||_{C^{\infty}_{F}}\, ,$$
where $K$ is a constant which depends only on the triplet and on the index $\eta_\infty$. Since the growth of the 
function does not change when the operator $\Lcal$ is applied, martingale property of $M^f$ can be done using Proposition 
\ref{prop:checkmg} with $\eta^*=\eta_\infty$. 

By iterating the above estimates, we get 
\begin{equation}\label{eq:est_iter_op}
||\Lcal^n f||_{C^\infty_{F}}\leq K^n ||f||_{C^{\infty}_{F}}\,.
\end{equation}
\endproof
\end{proposition}

The following result will be used in the following sections. 
\begin{corollary}\label{cor:rem}
Let $\Rcal_\nu f(y,s)$ be the remainder of order $\nu$ in \eqref{eq:taylorlevy}, 
$$\Rcal_{\nu}f(y,s):=\frac{s^{\nu+1}}{\nu!}\int_0^1(1-r)^\nu \expvBig{y}{\Lcal^{\nu+1}f(L_{rs})}dr\, .$$
Using the same notations of Proposition \ref{prop:levydynk}, for $s$ small, there exists a constant $C_{\nu,f}$ such that 
$$|\Rcal_\nu f(y,s)|\leq  s^{\nu+1}C_{\nu,f}F_{\eta^*}(y)\, .$$
\proof
Since, by assumption, $$|\partial^\alpha_x f(x)|\leq ||f||_{C^\infty_F} F_{\eta_\infty}(x),\mbox{ for all }\alpha\in\N^d\, ,$$ \eqref{eq:est_iter_op}
 implies
 \begin{equation}\label{eq:i}
 \expvBig{y}{|\Lcal^{\nu+1}f(L_{s})|}\leq K^{\nu+1}||f||_{C^\infty_F}\expvBig{y}{F_{\eta^*}(L_s)}.
 \end{equation}
 Moreover, using Dynkin's formula,
 \begin{align*}
 \expvBig{y}{F_{\eta_\infty}(L_s)}&=F_{\eta_\infty}(y)+\expvBig{y}{\int_0^s \Lcal F_{\eta_\infty}(L_u)du}\\
 &\leq F_{\eta_\infty}(y)+K\expvBig{y}{\int_0^s F_{\eta_\infty}(L_u)du}.
 \end{align*}
 By Gronwall's inequality
  \begin{align}\label{eq:ii}
 \expvBig{y}{F_{\eta_\infty}(L_s)}&\leq e^{Ks}F_{\eta_\infty}(y).
 \end{align}
 Combining \eqref{eq:i} and \eqref{eq:ii}
\begin{equation} 
 \expvBig{y}{|\Lcal^{\nu+1}f(L_s)|}\leq K^{\nu+1}e^{Ks}||f||_{C^\infty_F}F_{\eta_\infty}(y).
 \end{equation}
 We can estimate
 \begin{align*}
 |\Rcal_\nu f(y,s)|&\leq\frac{s^{\nu+1}}{\nu ! }\int_0^1 (1-r)^\nu\expvBig{y}{|\Lcal^{\nu+1} f(L_{sr})|}dr\\
 &\leq \frac{s^{\nu+1}}{\nu ! }K^{\nu+1}||f||_{C^\infty_F}F_{\eta_\infty}(y)\int_0^1 (1-r)^\nu e^{Krs}dr\\
 &\leq s^{\nu+1} C_{\nu,f}F_{\eta_\infty}(y)\\
  &\leq s^{\nu+1} C_{\nu,f}F_{\eta^*}(y).
 \end{align*}
\endproof
\end{corollary}

\section{Regularity results for Affine--type operators}
In the previous sections we made essentially two big assumptions on the function space $\Mcal$. 
The first one is continuity at time $t=0$ for the transition semigroup $t\mapsto P_tf(x)$ when $f\in \Mcal$ and the second one is martingale property of the process 
$M^f$ for $f\in\Mcal.$ We have seen that continuity of the transition semigroup can be achieved once we have enough integrability of the distributions so that a dominate convergence type theorem can be applied. Here we will restrict ourselves 
on subsets of weighted spaces where martingale property of $M^f$ holds.

\subsection{Affine Processes}
As a result of the previous section, we found that, under mild conditions on the function space $\Mcal$ and some continuity of the 
transition semigroup, we can derive differentiability in time of the function $u(t,x)$ via an iterated versions of the Dynkin's lemma. 
Since we have found a successful way to approach regularity in time, it is desirable to apply results in the previous section also for the analysis of regularity in space. 
This means that we seek for Markov processes whose transition semigroup $P_t f(x)$ makes sense also when seen as a function of the space parameter $x.$ 
In this section we focus on a class of Markov processes for which the transition semigroup as a function of $x$ can be identified with another Markov semigroup. 
We need to start with some additional notation

Define\begin{align*}
\Ccal&:=\{\eta:\Ucal\to\C\mbox{ of L\'evy--Khintchine form }\eqref{eq:infdiv}\},\\
\Ccal^*&:=\{\Psi:\Ucal\to\C^d\st \pi_I\Psi\in\Ccal^m\mbox{ and }\pi_J\Psi=A\pi_Ju\\
&\qquad\qquad\mbox{for some }A\in\R^{n\times n}\}.\nonumber
\end{align*}

\begin{theorem}[see \cite{2003}, \cite{nico_thesis}]\label{teo:bo}
Let $(\Psi_t)_{t\geq0}$ be a sequence of functions in $\Ccal^*$ which is differentiable in time satisfying 
\begin{equation}\label{eq:semiflow}
 \mbox{for any}\; s,t>0,\;\Psi_{t+s}(u)=\Psi_t(\Psi_s(u))
\end{equation}
and consider $d$ independent L\'evy processes 
$L^{\Psi_i},\;i=1,\ldots,d$ defined as
 in \eqref{eq:deflevy}. 
 
 \begin{enumerate}[(i)]
 \item For $f\in C_b(D)$ the family
\begin{equation}\label{eq:defsemi}
P_{\Psi}f(x):=\expvBig{}{f(L^{\Psi^{(1)}}_{x_1}+\ldots+L^{\Psi^{(d)}}_{x_d})},\qquad t\geq0.
\end{equation}
defines a stochastically continuous Markov process with state space $D$ with
\begin{equation}\label{eq:AP}
\expvBig{x}{e^{\la u,X_t\ra}}=e^{\la x,\Psi_t(u)\ra},\qquad t\geq0, u\in\Ucal.
\end{equation}
\item Let $X$ be a stochastically continuous Markov process with space space $D$ satisfying \eqref{eq:AP}. Then, for each $x\in D$ the process $X$ is a $\Pa^x$-semimartingale with semimartingale characteristics
\begin{align*}
A_t&=\int_0^t \alpha(X_{s^-})ds\\
B_t&=\int_0^t \beta(X_{s^-})ds\\
K(\omega,dt,d\xi)&=K(X_{t^-}(\omega,d\xi))dt
\end{align*}
where $\alpha,\,\beta$ and $K(x,d\xi)$ are functions of the form
\begin{align*}
\alpha(x)&=x_1\alpha_1+\ldots+x_m\alpha_d\\
\beta(x)&=x_1\beta_1+\ldots+x_d\beta_d\\
K(x,d\xi)&=x_1M_1(d\xi)+\ldots+x_mM_d(d\xi)
\end{align*}
and, for each $k=1,\ldots,d$ $\beta_k\in D$ and, for 
$i=1,\ldots,m,$ $\alpha_i\in S^+_d$  and $M_i$ are L\'evy measures on $D.$ 
\end{enumerate}

\end{theorem}

\begin{remark}
 Since $x\in D$ for all $j\in J$ one has $x_j\in\R.$ However observe that the definition \eqref{eq:defsemi} 
is still well defined because, due to the assumption on $\Ccal^*,$
$L^{\Psi_j}$ is a deterministic process. 
\end{remark}

The class of processes defined in Theorem \ref{teo:bo} is a subclass of the affine processes.

\begin{definition}[see \cite{2003}]
An \emph{affine process} is a time-homogeneous Markov process with state space $(D,\Bcal(D))$
$$X=(\Omega,(X_t)_{t\geq0},(\Fcal_t)_{t\geq0},(p_t)_{t\geq0},(\Pa^x)_{x\in D}),$$
satisfying
\begin{description}
 \item[{\bfseries stochastic continuity:}]\label{prop:stoc}for every $t\geq0$ and $x\in D,$ $\lim_{s\to t}p_s(x,\cdot)=p_t(x,\cdot)$ weakly,
\item[{\bfseries affine property:}]\label{prop:affine}there exist functions $\phi:\Rp{}\times\Ucal\to\C$ and $\Psi:\Rp{}\times\Ucal\to \C^d$ such that
\begin{equation}\label{eq:affineprop}
\expvBig{x}{e^{\la u,X_t\ra}}=\int_D e^{\la u,\xi\ra}p_t(x,d\xi)=e^{\phi(t,u)+\la \Psi(t,u),x\ra}
\end{equation}
 for all $x\in D$ and $(t,u)\in\Rp{}\times\Ucal.$
\end{description}

\end{definition}

It is possible to show that, up to an enlargement of a state space, given an affine process there is no loss of generality in assuming that \eqref{eq:affineprop} is replaced 
by \eqref{eq:defsemi}. See \cite{nico_thesis} for the details.

We will also use the following different characterizations:

\begin{proposition}[Proposition 4.13 \cite{nico_thesis}]
For any fixed $t\geq0$ and $x\in D$, the following are equivalent:
\begin{enumerate}[(i)]
\item\label{cor:1} $X^x_t$ is the value at time $t$ of an affine process starting from $x$ with  
\begin{equation}
\expvBig{x}{e^{\la u,X_t\ra}}=e^{\la x,\Psi(t,u)\ra},\qquad t\geq0, u\in\Ucal\,,
\end{equation}
\item\label{cor:2} there exists a Levy process $L^{(t,x)}$ such that
$$\expvBig{x}{e^{\la u,X_t\ra}}=e^{\la x, \Psi(t,u)\ra}=\expvBig{}{e^{\langle L^{(t,x)}_1,u\rangle}},\quad u\in\Ucal\, ,$$
\item\label{cor:3} in distribution it holds
$$X^x_t\stackrel{d}{=}L^{(t,e_1)}_{x_1}+\ldots+L^{(t,e_d)}_{x_d}\, ,$$
where $e_k$ are the canonical coordinates in $\R^d$ $$e_1:=(1,0,\ldots,0)^\top,\ldots,e_d:=(0,\ldots,0,1)^\top\, .$$
\end{enumerate} 
\end{proposition}

Each $L^{(t,e_k)},\;k=1,\ldots,d$, is a semimartingale with state space $D$,  by construction. Its semimartingale characteristics, 
relative to a truncation function $h$,  admit a version of the following form 
$$(s\mathtt b_k(t),s\sigma_k(t),s\nu_k(t,d\xi))\, ,$$ where 
$(\mathtt b_k(t),\sigma_k(t),\nu_k(t))$ is a L\'evy triplet for each $k=1,\ldots,d$.

\begin{notation}
Let $L^{(t,e_1)},\ldots,L^{(t,e_d)}$ be $d$ independent L\'evy processes each of 
them representing $X^{e_i}$ for $i=1,\ldots,d$. Henceforth, the following notation will be used:
\begin{itemize}
\item the extended generator of $L^{(t,e_i)}$ is denoted by $\Lcal^{(t,e_i)},$
\item its Markov semigroup is denoted by 
\begin{equation}\label{eq:evQ}
\expvBig{y}{f(L^{(t,x)}_s)}=v^{(t,x)}(s,y)=Q^{(t,x)}_sf(y),
\end{equation}
for all $ (t,x),(s,y)\in\Rp{}\times D$ and $f\in C_b(D).$
\end{itemize}
Unless differently specified, the notation $Q$ (resp. $\Lcal$) denotes the 
Markov semigroup (resp. the extended generator) of a L\'evy process, while the 
notations $P$ and $\Acal$ are reserved for the same quantities for an affine process.
\end{notation}

\begin{lemma}
Fix $t\geq0$ and $x\in D$. Suppose that  $\expvBig{x}{e^{\la y,X_t\ra}}<\infty$ for all $y\in\R^d$. Let $F_{\eta^*}$ be a weight functions such that, for some $y\in\R^d$, it holds
$$\sup_{x\in D} F_{\eta^*}(x)e^{-\la y,x\ra}<\infty.$$
Then, for all $f\in C^\infty_F$ with $\eta_\infty\leq\eta^*$ it holds $Q^{(t,x)}_s|f|<\infty$ for any fixed $(t,x)\in\Rp{}\times D$ and all $s\geq0.$
\proof
Let $q_t(x,\cdot)$ denote the distribution at time 1 of the L\'evy process $L^{(t,x)}$. Then, by assumption, it holds
$\int_D e^{\la y,\xi\ra}q_t(x,d\xi)<\infty$ for all $y\in\R^d$. Since $q_t(x,d\xi)$ is infinitely divisible, from Theorem 25.3 in \cite{sato_levy_1999}, we conclude that $L^{(t,x)}$ is a L\'evy process 
with L\'evy measure $\nu(t,x,\cdot)$ satisfying $$\int_{\{|\xi|\geq1\}} e^{\la y,\xi\ra}\nu(t,x,d\xi)<\infty\, ,$$
for all $y\in\R^d$. Moreover, finiteness of the exponential moments is not a time dependent property 
and therefore $\int_D e^{\la y,\xi\ra}q_t(sx,d\xi)<\infty$ for all $s\geq0$ and $y\in\R^d$. Hence, if $f$ satisfies the above mentioned growth condition,
\begin{align*}
 Q^{(t,x)}_s |f|(z)&=\int |f(z+\xi)|q_t(sx,d\xi)\\
&\leq \int |f(\xi)|q_t(sx,d\xi-z)\\
&\leq \int F_{\eta^*}(\xi)q_t(sx,d\xi-z).
\end{align*}
From spacial homogeneity we conclude that this quantity is finite.
\endproof
\end{lemma}

\begin{lemma}\label{lem:regreprelevy}
 For any $f:D\mapsto \R$ such that $\expvBig{x}{|f(X_t)|}<\infty$ for all $t\geq0$ and $x\in D$,  the following representation holds
$$u(t,x+hy)=\expvBig{x}{v^{(t,y)}(h,X_t)},\qquad h>0,\; y\in D\, .$$
\proof
For any $x,y\in D$ and $t,h>0$ it holds
\begin{eqnarray*}
 \expvBig{x+hy}{f(X_t)}&=&\int_D f(\xi)p_t(x+hy,d\xi)\\
&=&\int_D f(\xi)(p_t(x,\cdot)*p_t(hy,\cdot))(d\xi)\\
&=&\int_D f(\xi)\int_D p_t(x,d\xi-\eta)p_t(hy,d\eta)\\
&=&\int_D\int_D f(\xi+\eta)p_t(hy,d\eta)p_t(x,d\xi)\\
&=&\int_D \expvBig{}{f(\xi+L^{(t,y)}_h)}p_t(x,d\xi)\\
&=&\expvBig{x}{v^{(t,y)}(h,X_t)}\,.
\end{eqnarray*}
\endproof
\end{lemma}

\subsection{Results on \texorpdfstring{$C^\infty_{\textrm{pol}}$}{TEXT}}\label{sec:polygrowth}

In the field of weak approximation of SDE it is essential to have conditions which 
guarantee that the convergence error obtained for a certain numerical scheme 
in a small time horizon can be ``well propagated" up to a fixed time horizon. Regularity of the Kolmogorov equation with small initial data with polynomial growth allows 
to control this error (see \cite{alfonsi_high_2010}, \cite{talay} for example).

We first define the following function spaces: 
\begin{definition}
A function $f\in C^k_{\textrm{pol}}$ if
\begin{itemize}
\item $f\in C^k$
\item for all $\alpha$ multi-index with $|\alpha|\leq k$, there exist constants $C_\alpha$ and $\eta_\alpha$ such that
$$|\partial^\alpha f(x)|\leq C_\alpha(1+|x|^{2\eta_\alpha}),\quad\mbox{ for all }x\in D\, .$$
\end{itemize}
\end{definition}

In case the function $f$ is smooth we can extend the previous definition by taking all 
the possible derivatives and define
\begin{align*}
C^\infty_{\textrm{pol}}(D)=&\Big\{f\in C^\infty(D),\,\mbox{ for all } \alpha\in \N^d\,\;\exists C_\alpha>0, \eta_\alpha\in \N\\
&\quad\mbox{such that }  |\partial^\alpha f(x)|\leq C_\alpha(1+|x|^{2\eta_\alpha})\mbox{ for all }x\in D\Big\}\,.
\end{align*}

The first step is to check the basic properties. 
\begin{lemma}\label{lem:continzero}
Under the assumption that there exists a $T>0$ such that 
\begin{equation}\label{eq:momentcond}
\expvBig{x}{e^{\la y,X_T\ra}}<\infty\mbox{ for all  }y\in\R^d
\end{equation}
\begin{enumerate}[(i)]
\item  $t\mapsto P_t f$ is continuous at $t=0$  for all $f\in C^\infty_{\textrm{pol}},$
\item $P_tC^\infty_{\textrm{pol}}\subseteq C_{\textrm{pol}}.$
\end{enumerate}
\proof
We first check that $P_t|f|<\infty$ for all $f\in C^\infty_{\textrm{pol}}$. By assumption, there exist $C>0$ and $\eta>0$ 
such that $|f(x)|\leq C(1+|x|^{2\eta})=:C\Pi_{2\eta}(x)$,  where 
$\Pi_{2\eta}(x)$ is a polynomial of order $2\eta$. From the estimates in Theorem 2.10 in \cite{polypro} we know that there exists a 
$K>0$ such that 
$$\expvBig{x}{\Pi_{2\eta}(X_t)}\leq Ce^{Kt}F_{2\eta}(x),\qquad t\geq0\, ,$$
with $F_{2\eta}(x):=\left(1+\sum_{i=1}^d x^{2\eta}_i\right)$. Since $F_{2\eta}(x)\leq \Pi_{2\eta}(x)$,  we 
get integrability of $P_t f$. This, together with Lemma \ref{lem:continzerogeneral}, implies continuity at $t=0$ of $t\mapsto P_tf$. 
From the previous inequality we do also get the polynomial growth. 
 Finally, continuity follows from stochastic continuity, once the test functions are weighted with the weight function $F_\eta(x):=(1+|x|^{2\eta+1}).$
\endproof
\end{lemma}

\begin{theorem}\label{teo:FKpol}
Let $X$ be an affine process such that, for all $y\in\R^d$ and $x\in D$,  $\expvBig{x}{e^{\la y,X_T\ra}}<\infty$ for some fixed $T>0$. Then it holds
\begin{enumerate}[(i)]
 \item $\Acal C^{\infty}_{\textrm{pol}}\subseteq C^{\infty}_{\textrm{pol}}\,,$
\item for any $f\in C^\infty_{\textrm{pol}}$,  $P_tf$ solves the Kolmogorov's equation 
\begin{align*}
 \partial_t u(t,x)&=\Acal u(t,x)\,,\\
u(0,x)&=f(x)\,,
\end{align*}
for $(t,x)\in[0,T]\times D,$
\item for any $f\in C^\infty_{\textrm{pol}}$ and $\nu\in\N$ the following expansion of 
the transition semigroup holds for $(t,x)\in[0,T]\times D:$
\begin{equation*}
 \expvBig{x}{f(X_t)}=f(x)+\sum_{k=1}^\nu \frac{t^k}{k!}\Acal^k f(x)+\Rcal_\nu f(x,t),
\end{equation*}
where $\Rcal_{\nu}f(x,t)$ is a remainder of order $\Oh(t^{\nu+1})$. 
\end{enumerate}
\proof
Using linearity in $x$ of the coefficients, we decompose
$$\Acal f(x)=\sum_{i=1}^d x_i \Acal^{(i)} f(x)\, ,$$
where each $\Acal^{(i)}$ is an operator of L\'evy-type. For every $i=1,\ldots,m$,
\begin{align*}
 \Acal^{(i)} f(x)&=\sum_{k=1}^d(\beta_i)_k\partial_{x_k} f(x)+\frac{1}{2}\sum_{k,h=1}^d(\alpha_i)_{kh}\partial^{2}_{x_k x_h} f(x)\\
&+\int\left( f(x+\xi)-f(x)-\la h(\xi),\nabla f(x)\ra\right)M_i(d\xi)\,.
\end{align*}
From the integrability assumption, it follows that $\int_{\{|\xi|\geq1\}}e^{\la y,\xi\ra}M_i(d\xi)<\infty$,  for all $i=1,\ldots,m$ and $y\in\R^d$ (see 
Theorem 2.14 in \cite{expmom}). Hence
\begin{align*}
&\left|\int_{\{|\xi|>1\}}\left(f(x+\xi)-f(x)\right)M_i(d\xi)\right|\leq C_1(1+|x|^{2\eta_1})\int_{\{|\xi|>1\}}|\xi|M_i(d\xi).\,\\
&\left|\int_{\{|\xi|\leq1\}}\left( f(x+\xi)-f(x)-\la \xi,\nabla f(x)\ra\right)M_i(d\xi)\right|\\
&\qquad\qquad\qquad\qquad\qquad\qquad\qquad\leq C_2(1+|x|^{2\eta_2})\int_{\{|\xi|\leq1\}}|\xi|^2M_i(d\xi)\,.&
\end{align*}
Moreover, since $f\in C^{\infty}_{\textrm{pol}}$ there exist two constants $C$ and $E$ such that 
\begin{equation}\label{def:CE}
|f(x)|+\sum_{i=1}^d|\partial_{x_i} f(x)|+\sum_{i,j=1}^d |\partial^2_{x_ix_j}f(x)|\leq C(1+|x|^{2E})\,.
\end{equation}

Combining the bound on the diffusive part and the jump part, we conclude that there exist two constants $K$ and $E$ such that 
$$|\Acal^{(i)} f(x)|\leq K(1+|x|^{2E}),\quad\mbox{for all }i=1,\ldots,d\, .$$
Then 
$$|\Acal f(x)|\leq K(1+|x|^{2E})|x|\leq\overline K(1+|x|^{2(E+1)})\, .$$
This concludes (i).

In order to apply Theorem \ref{teo:FK}, we still need to prove that, for 
$f\in C^\infty_{\textrm{pol}}$,
\begin{itemize}
\item  process $M^f$ is a true martingale, 
\item $t\mapsto P_t f(x)$ is continuous at $t=0$.
\end{itemize}
We start with the martingale property. Let $(C,E)$ be the couple defined in \eqref{def:CE}. Since $|\Acal f|\leq K (1+|x|^{2(E+1)})$,  in order to apply Proposition \ref{prop:checkmg}, it remains to check that 
$$\expvBig{x}{\sup_{t\in[0,T]}|X_t|^{2(E+1)}}<\infty\, .$$

By Lemma 2.17 in \cite{polypro}, there exist two constants $K$ and $C$ such that 
$$\expvBig{x}{\sup_{t\in[0,T]}|X_t|^{2(E+1)}}\leq Ke^{CT}\, ,$$
if the kernel $K(x,d\xi)=x_1M_1(d\xi)+\ldots+x_dM_d(d\xi)$ satisfies
\begin{equation}\label{eq:polycond}
\int_{\{|\xi|\geq1\}} |\xi|^{2(E+1)}K(X_t,d\xi)\leq (1+|X_t|^{2(E+1)}).
\end{equation} 
By assumption, for all $i=1,\ldots,d,$
$$\int_{\{|\xi|\geq1\}}|\xi|^{2(E+1)}M_i(d\xi)<\infty\, .$$
Hence 
$$\int_{\{|\xi|\geq1\}} |\xi|^{2(E+1)}K(X_t,d\xi)=\int_{\{|\xi|\geq1\}} |\xi|^{2(E+1)}\sum_{i=1}^d X^{(i)}_t M_i(d\xi)\leq C(1+|X_t|)$$
and therefore \eqref{eq:polycond} holds. Now that martingale property has been proved, it remains to check that $t\mapsto P_tf$ is continuous at $t=0$. This follows from the integral equation
$$P_t f(x)=f(x)+\int_0^t P_s \Acal f(x)ds\, ,$$
paired with the fact that $\Acal$ maps 
$C^\infty_{\textrm{pol}}$ into $C^\infty_{\textrm{pol}}$. From Theorem \ref{teo:FK} we conclude that for any $f\in C^\infty_{\textrm{pol}}$ and $P_tf$ solves the 
Kolmogorov's equation. Finally (iii) follows as in Proposition \ref{prop:dynkiniter} by considering that 
$\Acal$ maps functions in $C^{\infty}_{\textrm{pol}}$ into functions in $C^{\infty}_{\textrm{pol}}\,.$
\endproof
\end{theorem}

To show
$$P_t C^\infty_{\textrm{pol}}\subseteq C^\infty_{\textrm{pol}}$$

we will consider the decomposition 
$$X^{x+he_i}\stackrel{law}{=} X^x+\widetilde X^{he_i},\;h>0, i=1,\ldots,d\, ,$$
where $e_1,\ldots,e_d$ are the basis elements in $\R^d$ and $\widetilde X^{he_i}$ is an independent copy of the process $X$ starting from $he_i$.

\begin{theorem}\label{teo:takederonce}
Suppose that there exists a $T>0$ such that 
$\expvBig{x}{e^{\la y, X_T\ra}}<\infty$ for all $y\in\R^d$. Then, all the partial derivatives exist and are continuous. Moreover it holds
 \begin{equation}\label{eq:takederonce}
\partial_{e_i}u(t,x)=\expvBig{x}{\Lcal^{(t,e_i)}f(X_t)},\qquad t\geq0, x\in D.
\end{equation}

\proof
Let $p_t(x,\cdot)$ be the distribution of $X^x_t$ for $t\in[0,T]$. From Lemma 4.2 (c) in \cite{expmom}, for all $t\in[0,T],$
$$\int_D e^{\la y,\xi\ra}p_t (x,d\xi)<\infty\, .$$
Recall that $p_t(x,\cdot)$ are infinitely divisible measures. Denote by $\nu_t(x,\cdot)$ their L\'evy measures. From from Theorem 25.3 in \cite{sato_levy_1999}, it follows that 
$$\int_{\{|\xi|\geq1\}} e^{\la y,\xi\ra}\nu_t(x,d\xi)<\infty\,,$$ for all $y\in\R^d$ and $t\in[0,T]$. Due to linearity in $x$,  Proposition \ref{prop:levydynk} can be applied to each $L^{(t,e_i)}$. This allows us to get the following approximate for $h$ small
\begin{eqnarray*}
 u(t,x+he_i)-u(t,x)&=&\expvBig{x}{v^{(t,e_i)}(h,X_t)-v^{(t,e_i)}(0,X_t)}\\
&\leq&h\expvBig{x}{\Lcal^{(t,e_i)}f(X_t)}+\Oh(h^2).
\end{eqnarray*}
By taking the limit of the finite differences, we get existence of the derivatives. Additionally, since $\Lcal^{(t,e_i)}$ maps
 functions in $C^{\infty}_{\textrm{pol}}$ into functions in $C^{\infty}_\textrm{pol}$,  Lemma \ref{lem:continzero} 2. leads to continuity of the derivatives. 
\endproof
\end{theorem}

 Higher order partial derivatives can be taken analogously by applying \eqref{eq:takederonce}
 several times:
\begin{proposition}\label{prop:highder}
Let $\alpha$ be a multi--index and $f\in C^{\infty}_{\textrm{pol}}$. Then, under the same assumptions as in Theorem \ref{teo:takederonce}, 
\begin{equation}\label{eq:takeder}
 \partial^\alpha_x u(t,x)=\expvBig{x}{(\Lcal^{(t,e_1)})^{\alpha_1}\cdots(\Lcal^{(t,e_d)})^{\alpha_d}f(X_t)}, 
\end{equation}
and therefore all the derivative exists. Moreover, they are continuous in $x$.
\proof
The representation of the partial derivatives follows analogously to the case when $|\alpha|=1$. 
From \eqref{eq:takeder} we can write 
$\partial^\alpha_x u(t,x)=P_t g(x)$,  with  
$g(x):=(\Lcal^{(t,e_1)})^{\alpha_1}\cdots(\Lcal^{(t,e_d)})^{\alpha_d} f(x)\in C^\infty_{\textrm{pol}}$. The result holds since $P_t C^\infty_{\textrm{pol}}\subseteq C_\textrm{pol}.$
\endproof
\end{proposition}

\begin{proposition}
 Let $X^x$ be an affine process satisfying the assumptions in Theorem \ref{teo:takederonce}. Then 
given a function $f\in C^\infty_{\textrm{pol}}$,  
$P_t f(x)$ is again in $C^\infty_{\textrm{pol}}$ for all $t\geq0.$
\end{proposition}

\subsection{Results on \texorpdfstring{$C^\infty_{\textrm{exp}}$}{TEXT}}\label{sec:expgrowth}

In some cases, see next section, polynomial growth is not enough and we need exponential growth of the test function.

In this section we will work with the following function spaces: 
\begin{definition}
A function $f\in C^k_{\textrm{exp}}$ if
\begin{itemize}
\item $f\in C^k$
\item for all $\alpha$ multi-index with $|\alpha|\leq k$, there exist constants $C_\alpha$, $\theta_\alpha>0$ such that
$$|\partial^\alpha f(x)|\leq C_\alpha\cosh(\theta_\alpha|x|),\quad\mbox{ for all }x\in D\, .$$
\end{itemize}
\end{definition}

In case the function $f$ is smooth we can extend the previous definition by taking all 
the possible derivatives and define
\begin{align*}
C^\infty_{\textrm{exp}}(D)=&\Big\{f\in C^\infty(D),\,\mbox{ for all } \alpha\in \N^d\,\;\exists C_\alpha>0, \theta_\alpha>0\\
&\quad\mbox{such that }  |\partial^\alpha f(x)|\leq C_\alpha\cosh(\theta_\alpha|x|)\mbox{ for all }x\in D\Big\}\,.
\end{align*}

The analogous of Lemma \ref{lem:continzero} holds with a slight modification in the proof
\begin{lemma}
Under the exponential moment condition \eqref{eq:momentcond} 
\begin{enumerate}[(i)]
\item  $t\mapsto P_t f$ is continuous at $t=0$  for all $f\in C^\infty_{\textrm{exp}},$
\item $P_tC^\infty_{\textrm{exp}}\subseteq C_{\textrm{exp}}.$
\end{enumerate}
\proof
We first check that $P_t|f|<\infty$ for all $t\geq0$. In the following the constant $C$ may be different from line to line. 
By assumption there exists $C,\theta>0$ such that
$$\expvBig{x}{|f(X_t)|}\leq \expvBig{x}{C\cosh(\theta|X_t|)}\leq C\expvBig{x}{e^{\theta |X_t|}}\,.$$
Writing the last term componentwise we see that 
$$\expvBig{x}{e^{\theta X_t}}\leq \expvBig{x}{\exp\left({\theta\sum_{i=1}^d |X^{(i)}_t|}\right)}\leq\expvBig{x}{\prod_{i=1}^d\left(e^{\theta X^{(i)}_t}+e^{-\theta X^{(i)}_t}\right)}$$
which can be written as a sum of finite terms of the form $\expvBig{x}{e^{\ell(X_t)}}$, where $\ell$ is a linear function in $\R^d$. 
Due to the exponential moment condition, we conclude, using Lemma \ref{lem:continzerogeneral}, the continuity at $t=0$ of $t\mapsto P_tf$. 
The exponential growth follows applying the results in \cite{expmom}. For each $y\in\R^d$, Theorem 2.14 in \cite{expmom} gives the existence of a $C^1$ function 
$q:t\mapsto q(t,y)$ from $[0,T]$ to $\R^d$ such that
$$\expvBig{x}{e^{\la u,X_t\ra}}=e^{\la x,q(t,y)\ra}$$
holds for all $(t,x)\in[0,T]\times D$. Therefore, using the same argument as before
$$\expvBig{x}{f(X_t)}\leq C\expvBig{x}{e^{\ell(X_t)}}=e^{\la x,q(t,y)\ra}\,,$$
for some $y\in\R^d$, from where the exponential growth. Finally, continuity goes as in the polynomial growth case. 
\endproof
\end{lemma}

\begin{theorem}\label{teo:FKexp}
Let $X$ be an affine process such that, for all $y\in\R^d$ and $x\in D$,  $\expvBig{x}{e^{\la y,X_T\ra}}<\infty$ for some fixed $T>0$. Then it holds
\begin{enumerate}[(i)]
 \item $\Acal C^{\infty}_{\textrm{exp}}\subseteq C^{\infty}_{\textrm{exp}}\,,$
\item for any $f\in C^\infty_{\textrm{exp}}$,  $P_tf$ solves the Kolmogorov's equation 
\begin{align*}
 \partial_t u(t,x)&=\Acal u(t,x)\,,\\
u(0,x)&=f(x)\,,
\end{align*}
for $(t,x)\in[0,T]\times D,$
\item for any $f\in C^\infty_{\textrm{exp}}$ and $\nu\in\N$ the following expansion of 
the transition semigroup holds for $(t,x)\in[0,T]\times D:$
\begin{equation*}
 \expvBig{x}{f(X_t)}=f(x)+\sum_{k=1}^\nu \frac{t^k}{k!}\Acal^k f(x)+\Rcal_\nu f(x,t),
\end{equation*}
where $\Rcal_{\nu}f(x,t)$ is a remainder of order $\Oh(t^{\nu+1})$. 
\end{enumerate}
\proof The proof goes as in Theorem \ref{teo:FKpol}. Indeed, 
\begin{align*}
&\left|\int_{\{|\xi|>1\}}\left(f(x+\xi)-f(x)\right)M_i(d\xi)\right|\leq C_1\cosh(\theta_1|x|)\int_{\{|\xi|>1\}}|\xi|M_i(d\xi).\,\\
&\left|\int_{\{|\xi|\leq1\}}\left( f(x+\xi)-f(x)-\la \xi,\nabla f(x)\ra\right)M_i(d\xi)\right|\\
&\qquad\qquad\qquad\qquad\qquad\qquad\qquad\leq C_2\cosh(\theta_2|x|)\int_{\{|\xi|\leq1\}}|\xi|^2M_i(d\xi)\,.&
\end{align*}
Moreover, since $f\in C^{\infty}_{\textrm{exp}}$ there exist two constants $C$ and $\Theta$ such that 
\begin{equation}
|f(x)|+\sum_{i=1}^d|\partial_{x_i} f(x)|+\sum_{i,j=1}^d |\partial^2_{x_ix_j}f(x)|\leq C\cosh(\Theta |x|)\,.
\end{equation}
Since for all $\theta>0$ there exists a constant $C>0$ such that $|x|\cosh(\Theta|x|)\leq C\cosh(2\Theta|x|)$, combining the bound on the diffusive part and the jump part, 
we conclude (i). In order to get martingale property of $M^f$ with $f\in C^\infty_{\textrm{exp}}$ we can use Proposition \ref{prop:checkmg}.

Finally (iii) follows as in Proposition \ref{prop:dynkiniter} by considering that 
$\Acal$ maps functions in $C^{\infty}_{\textrm{exp}}$ with good couple $(C,\theta)$ into functions in $C^{\infty}_{\textrm{exp}}$ with good couple $(C,2\theta)$.
\endproof
\end{theorem}

In analogy with the previous section we obtain the following result

\begin{proposition}
Given a function $f\in C^\infty_{\textrm{exp}}$,  
$P_t f(x)$ is again in $C^\infty_{\textrm{exp}}$ for all $t\geq0.$
\end{proposition}

\section{Applications}
\subsection{Analysis of convergence rates}
Suppose that $X$ is an affine process and $H$ is a functional acting on the paths. In mathematical finance, we may interpret $H$ as a payoff function, 
possibly depending on the whole path up to a fix time 
$T$, for a contract written on $X$. In order to price this constract, we need the compute the quantity $\expvBig{x}{H(X_t,\;t\in[0,T])}$, where $\expvBig{x}{}$ is the expectation taken 
under the pricing measure. In only few circumstances it is possible 
to derive a closed formula for this quantity and therefore it would be helpful to approximate the problem by taking an approximation of the path $X$ on a fixed time partition and 
then use Monte Carlo. 
Suppose that $\{t_0,\ldots,t_N\}$ is a uniform partition of $[0,T]$ with mesh size $h>0$. We want to find an approximating sequence $\{\tc X^x_{t_k}\}_{l=0,\ldots,N}$ such that
\begin{itemize}
\item[-] $\tc X^x_{t_0}=x$,
\item[-] it is a $\nu$-order approximation of $X^x$ in the sense that, for every $f$ smooth function, there exists a constant $C(f,x)$ and an index $\nu$ such that 
\begin{equation}\label{eq:globalerror}
\big|\E^{x}\left[f(X_T)\right]-\E\left[ f(\tc X^x_{t_N})\right]\big|\,\leq\, C(f,x) h^\nu.
\end{equation}
\end{itemize}
A possible way to define the approximating sequence is the following: define $\tc X_h$ as a random variable such that 
\begin{equation}\label{eq:localerror}
\Big|\expvBig{x}{f(X_{h})}-\expvBig{}{f(\tc X^x_{h})}\Big|\leq h^{\nu+1} C(f,x),\quad h\in (0,h_0),\;\mbox{with }h_0>0\,.
\end{equation}
Then, under a set of conditions we are going to specify soon, if $\{(\tc X^x_{h})^{k}\}_{k=1,\ldots,N}$ are independent copies of $\tc X^x_{h}$ the 
following iterative procedure
\begin{align*}
 \tc X^x_{t_0}&=x,\\
\tc X^x_{t_1}&=(\tc X^x_h)^1,\\
\tc X^x_{t_2}&={(\tc X^y_h)^2}_{\Big|y=\tc X^x_{t_1}},\\
\vdots&\\
\tc X^x_{t_N}&={(\tc X^y_h)^N}_{\Big|y=\tc X^x_{t_{N-1}}},\\
\end{align*}
defines a $\nu$-order scheme for $X^x$. To be more precise, identify $\expvBig{x}{f(X_t)}$ with the semigroup associated with the affine process $X^x$ 
acting on the function $f$ and define
$$Q_hf(x)=\expvBig{}{f(\tc X^x_h)}.$$
Then, by the previous iterative construction of the approximating sequence, it follows that
$$\expvBig{}{f(\tc X^x_{t_n})}=Q^n_h f(x),\qquad n=1,\ldots,N,$$
where $Q^n_h$ is the operator obtained by taking the $n$th--composition of $Q_h$. 
Then, the estimate in \eqref{eq:localerror} can be transported into \eqref{eq:globalerror} by means of the following telescopic sum
$$P_T f(x)-Q^N_h f(x)=\sum_{k=1}^N Q_h^{N-k}(Q_h-P_h)P_{t_k}f(x),$$
as long as all the three terms in the sum are well defined and well propagate the local error.

We first require that, small time horizon, 
$P_h f$ and $Q_h f$ admit the same expansion of some order $\Oh(h^{\nu+1})$. A second property is a regularity condition of the semigroup which implies preservation of the local error. 
The implication of these two requirements is that the function $f$ lies on a space $\Hcal$ which satisfies
\begin{description}
\item[H1)] for all $f\in\Hcal,$ there exists a $\nu\in\N$ such that $(P_h-Q_h)f\simeq\Oh(h^{\nu+1})$, for all $h\in(0,h_0)$,
\item[H2)] for all $f\in\Hcal$, $P_tf\in\Hcal$ and $Q_tf\in\Hcal.$ 
\end{description}
While the first property concerns the discretization scheme, the second one relates to the application of the semigroup on $\Hcal$. For weak approximation schemes of Markov processes, 
the space $\Hcal$ is given by $C^\infty_{\textrm{pol}}$. Here, for example, we state the general guideline from \cite{alfonsi_high_2010}. 

\begin{theorem}[see Theorem 1.9 in \cite{alfonsi_high_2010}]\label{teo:alfonsi} Let $X$ be an affine process with $X_0=x$ satisfying the exponential moment condition \eqref{eq:momentcond}. Denote by $\tc X^x$ its approximation on the uniform partition 
$\{t_0,\ldots,t_N\}$. Assume that
\begin{enumerate}[i)]
\item for all $h\in(0,h_0)$ and $\alpha\in\N$ there exists a constant $C_\alpha$ such that 
$$\expvBig{x}{|\tc X^x_h|^\alpha}\leq |x|^\alpha(1+C_\alpha)h+C_\alpha h\,,$$
\item for all $h\in(0,h_0)$ and $f\in C^\infty_{\textrm{pol}}$ there exist two constants $C,E>0$ such that
$$|P_h f(x)-Q_h f(x)|\leq Ch^{\nu+1}(1+|x|^E)\,,$$
\item for all $f\in C^\infty_{\textrm{pol}}$ the function $u(t,x):=P_t f(x)$ is well defined for $(t,x)\in [0,T]\times D$, is a smooth solution of 
$\partial_t u(t,x)=\Acal u(t,x)$ such that, for all $\alpha\in\N^{d+1}$ multi--index it holds
\begin{equation}
\mbox{for all }(t,x)\in[0,T]\times D,\quad |\partial^\alpha_{(t,x)}u(t,x)|\leq K_\alpha(T)(1+|x|^{\eta_\alpha(T)})\, , 
\end{equation}
where $K_\alpha(T)$ and $\eta_\alpha(T)$ are positive constants depending on the time horizon $T$ and the 
order of derivative $\alpha.$
\end{enumerate}
Then, there exists $K(x,T)>0$ such that
$$|P_T f(x)-Q_h^N f(x)|\leq K(x,T) h^\nu\,.$$ 
\end{theorem}

\begin{remark}
 The above result is stated in \cite{alfonsi_high_2010} only for It\^o process with no jumps and with coefficients in $C^\infty_{\textrm{pol}}$. However, the skeleton of the proof 
is given by the previous consideration on the 
expansion of the error via telescopic sum. The same steps in \cite{alfonsi_high_2010} can be carried also in the case of general affine processes using the results from Section \ref{sec:polygrowth}.
\end{remark}

In the above theorem, the first two conditions are specific properties of the approximation scheme which is chosen. Observe that, for affine processes, a standard Euler--Maruyama scheme can run into troubles due to the 
square--root term appearing on the volatility coefficient. High order schemes, like the one proposed in \cite{alfonsi_high_2010} can be applied. Anyway, suppose we are given an approximation which satisfies 
conditions (i) and (ii) in the above result. Condition (iii) is an intrinsic property of the stochastic process to approximate. When $X$ is an affine process, we can use the results from Section \ref{sec:polygrowth}.

\begin{theorem}\label{teo:finalteo}
 Let $f\in C^\infty_{\textrm{pol}}$. Then, the function $u:\Rp{}\times D\to\R$ defined by 
$u(t,x)=\expvBig{x}{f(X_t)}$ is smooth, with all derivatives satisfying the following property:
\begin{equation}
\mbox{for all }(t,x)\in[0,T]\times D,\quad |\partial^\alpha_{(t,x)}u(t,x)|\leq K_\alpha(T)(1+|x|^{2\eta_\alpha(T)})\, , 
\end{equation}
where $K_\alpha(T)$ and $\eta_\alpha(T)$ are positive constants depending only on the time horizon $T$ and the 
order of derivative $\alpha.$
\proof
Let $\alpha$ be a multi--index running over the mixed derivatives with respect to time and space. Split $\alpha=(\alpha_0,\overline\alpha)$ where $\alpha_0\in \N$ 
is the order of derivation in time and $\overline\alpha\in\N^d$ is the multi--index for the derivatives in space. Clearly 
$\alpha_0=|\alpha|-|\overline\alpha|$. By induction on $\alpha_0$,  when $\alpha_0=1$,  
$$\partial_t \partial^{\overline\alpha}_x u(t,x)=\Acal\partial^{\overline\alpha}_x u(t,x)\, .$$
We need to check that the bound $\partial^{\overline\alpha}_x u(t,x)$ can be taken uniformly for $t\in[0,T].$
From the representation \eqref{eq:takeder} and Corollary \ref{cor:rem},
\begin{eqnarray*}
 \partial^{\overline\alpha}_x u(t,x)&=&\expvBig{x}{(\Lcal^{(t,e_1)})^{\alpha_1}\cdots(\Lcal^{(t,e_d)})^{\alpha_d}f(X_t)}\\
&\leq& \prod_{k=1}^d K(t,x_k)^{\alpha_k}C(1+|x|^{2\eta_\alpha})\, ,
\end{eqnarray*}
where $K(t,x_k),\;k=1,\ldots,d$,  depends only on the L\'evy triplet of the distribution of  $X^{x_k}_t$. Due to linearity in $x$,  for all $t\in[0,T]$ there exists a constant $K$ depending on the time horizon $T$ such that 
$K(t,x_k)\leq K(T)x_k$. Therefore
$$\prod_{k=1}^d K(t,x_k)^{\alpha_k}\leq\prod_{k=1}^d K(T)^{\alpha_k}x_k^{\alpha_k}\leq K_{\alpha}(T)|x|^{|\overline\alpha|}\, .$$
Using the fact that, for all $u>0$ and $p,m\in\N$,
$$u^p(1+u^{2m})\leq (1+u^{2(m+p)})$$ we get
\begin{eqnarray*}
 \partial^{\overline\alpha}_x u(t,x)&\leq& K_{\alpha}(T)(1+|x|^{2(\eta_\alpha+|\overline\alpha|)})\, ,
\end{eqnarray*}
and the result follows by Theorem \ref{teo:FKpol} (i). Now suppose that, for all $\alpha_0=1,\ldots,|\alpha|-|\overline\alpha|-1$ it holds
$$|\partial^{\alpha_0}_t\partial^{\overline\alpha}_x u(t,x)|\leq K_\alpha (T)(1+|x|^{2\eta_\alpha(T)}),\;\mbox{for all }t\in[0,T]\, .$$
Then 
\begin{eqnarray*}
 \partial^{\alpha_0+1}_t\partial^{\overline\alpha}_x u(t,x)&=&\partial_t\left(\partial^{\alpha_0}_t\partial^{\overline\alpha}_x u(t,x)\right)\\
&=&\Acal \partial^{\alpha_0}_t\partial^{\overline\alpha}_x u(t,x)
\end{eqnarray*}
and, by the induction step, the result follows again as an application of Theorem \ref{teo:FKpol}.
\endproof

\end{theorem}

\begin{remark}For numerical applications, the existence of some order moments may be sufficient. The existence of exponential moments in 
the previous section has been used mainly to guarantee that moments of any order are well defined and finite. In practice, the approximation $Q_h f$ is defined by means of 
a truncation of the series in 
Theorem \ref{teo:FKpol} (iii) to a fixed order $\nu$. Recall that the expansion makes sense as long as $f$ is a function such that $f\in\Dcal(\Acal^{k})$ and 
$\Acal^k f\in\Mcal$ for all $k=1,\ldots,\nu+1$ (compare with Proposition \ref{prop:dynkiniter}). In turn, the martingale property can be check using Proposition \ref{prop:checkmg} 
requiring the existence of as many moments as necessary. Just as an example, consider the following situation. Let $f\in C^\infty_\textrm{pol}$ and let $(C_\alpha,\eta_\alpha)_{\alpha\in\N^d}$ be 
its good sequence. Introduce the notation
$$\overline\eta_k(f):=\max_{|\alpha|\leq k}\eta_\alpha\,.$$ 
Under the assumption that 
$$\int_{\{|\xi|>1\}}|\xi|^{2\overline \alpha_1}M_i(d\xi)<\infty,\quad\mbox{for all}\quad i=1,\ldots,d\,,$$
we know that $|P_t f(x)|\leq K(1+|x|^{2\overline \alpha_1})$ for all $t\in[0,T]$. This follows from \cite{polypro}. 
At the same time, if $g\in C^\infty_\textrm{pol}$ is such that $g,\Acal g,\Acal^2 g\in\Mcal$, then 
$$Q_h(x):=g(x)+\Acal g(x)h$$
defines a local second order scheme for $X$ (compare with Theorem \ref{teo:FKpol}). In terms of integrability conditions, the above assumption translates into 
the existence of moments up to order $4(\overline\eta_6(g)+1)$. Since the function $g(x):=P_t f(x)$ is in $C^\infty_{\textrm{pol}}$ with 
$\overline\eta_6(g)\leq\overline\eta_6(f)+6$, we conclude that, as long we as require the existence of $4(\overline\eta_6(f)+7)$ moments, 
$$(Q_h -P_h)P_{t_k} f(x)=\Rcal_2 P_{t_k} f(x),$$
where $\Rcal_2 P_{t_k} f(x)$ is a reminder of order $\Oh(h^2)$.
Finally, if $\tc X^x_h$ satisfies the property (i) in Theorem \ref{teo:alfonsi} with $\alpha\leq 2(\overline\eta_6(f)+6)$ then 
$$\expvBig{}{|\tc X^x_{t_{n}}|^\alpha}\leq  \expvBig{}{|\tc X^x_{t_{n-1}}|^\alpha}(1+C_\alpha)h+C_\alpha h$$
from where we conclude, from an application of the discrete version of Grownwall's lemma, 
$$\expvBig{}{|\tc X^x_{t_{n}}|^\alpha}\leq e^{KT}|x|^\alpha,\quad \mbox{for all } n=1,\ldots,N\,.$$
Therefore, plugging these estimates in the telescopic sum, we get
$$|P_T f(x)-Q^N_h f(x)|\leq \sum_{k=1}^N |Q^{N-k}_h \Rcal_2 f(x)|\leq h e^{KT}|x|^{2(\overline\eta_6(f)+6)}\,.$$
\end{remark}





\bibliographystyle{alpha}
\bibliography{sample.bib}







\end{document}